\documentclass[a4paper, twoside, 11pt, USenglish]{amsart}

	\usepackage[utf8]{inputenc}
	\usepackage[T1]{fontenc}
	\usepackage{multicol}

	\usepackage{amsmath,amssymb}
	\usepackage{mathtools}
	\usepackage[all]{xy}
	\usepackage{tikz-cd}

	\usepackage[textwidth=14.5cm,textheight=21.5cm,heightrounded,hcentering]{geometry}
	\usepackage[lighttt]{lmodern} 
	\usepackage[protrusion=true,expansion=true,kerning, spacing]{microtype}
		\microtypecontext{spacing=nonfrench}
	\usepackage{csquotes}

\usepackage{graphicx}
\usepackage{xcolor}
	\definecolor{Mycolor1}{HTML}{83c5be}
	\definecolor{Mycolor2}{HTML}{fa9483}
	\definecolor{Mycolor3}{HTML}{2d4057}
	\definecolor{Mycolor4}{HTML}{4097aa}
	\definecolor{e-mail}{rgb}{0,.40,.80}
	\definecolor{reference}{rgb}{.20,.60,.22}
	\definecolor{mrnumber}{rgb}{.80,.40,0}
	\definecolor{citation}{rgb}{0,.40,.80}
	\definecolor{gris25}{gray}{0.45}
	\definecolor{MyBlue}{rgb}{0,.12,.35}
\usepackage{wrapfig}
\usepackage{subcaption}
\usetikzlibrary{braids}

\usepackage{algpseudocode}
\usepackage{algorithm}
	
	\usepackage{calc} 
	\usepackage{caption}
	\DeclareMathOperator{\Gal}{Gal}
	\newcommand{\Q}{\mathbb{Q}}
	\newcommand{\Gq}{\Gal(\bar{\Q}/\Q)}
	\DeclareMathOperator{\Aut}{Aut}

\usepackage[USenglish]{babel}

	\usepackage{ellipsis}

	\author{Benjamin \textsc{Collas}}
	\author{\textsc{Murotani} Takahiro}
	\author{\textsc{Yamaguchi} Naganori}
		
	\title{Symmetries of spaces and numbers -- anabelian geometry}

	\usepackage[explicit]{titlesec}
		\titleformat{\section}[block]{\filcenter\scshape}{\thesection.\,#1}{1em}{}[]
		\titleformat{\subsection}[runin]{\bfseries}{\thesubsection.\,#1}{0em}{}[.\quad]
		\titleformat{name=\subsubsection}[runin]{\itshape}{\S\,#1}{0em}{}[.\quad]
	
	\usepackage{titletoc}
		\titlecontents{section}[0em]{\addvspace{.7pc}\bfseries\fontfamily{phv}\selectfont\scshape}{}{}{}[]
		\titlecontents{subsection}[1.5em]{\addvspace{0pc}\bfseries}{\thecontentslabel\ -\ }{}{\titlerule*[.8pc]{.}\contentspage}[]

		\titlecontents*{subsubsection}[1.5em]{}{\contentslabel[\relax]{0em}}{}{}[\,\textbullet\,][.]

\usepackage{filecontents}

\usepackage{csquotes}
\usepackage[
bibencoding=utf8,%
giveninits=true,
bibstyle=alphabetic,%
citestyle=alphabetic,%
ibidtracker=true,
sorting=nyt,
language=english,%
isbn=false,%
doi=false,%
eprint=true,
hyperref=true,
url=false,
minbibnames=4,
maxbibnames=4,
backend=biber,
datamodel=ext-eprint 
]{biblatex}
\bibliography{Anabelian-Snapshot}

\AtBeginBibliography{\small\vspace{-.5em}}
\defbibenvironment{bibliography}
{\list
	{\bfseries [\printfield[labelnumberwidth]{labelalpha}\printfield[labelnumberwidth]{extraalpha}]}
	{\setlength{\labelwidth}{0pt} 
		\setlength{\leftmargin}{0pt} 
		\setlength{\labelsep}{4pt}
		\setlength{\itemsep}{4pt}
		\setlength{\parsep}{2pt}} 
	}
{\endlist}
{\item}
\newbibmacro{string+doiurl}[1]{%
	\iffieldundef{doi}{%
		\iffieldundef{url}{%
			#1%
		}{%
			\href{\thefield{url}}{#1}%
		}%
	}{%
		\href{https://doi.org/\thefield{doi}}{#1}%
	}%
}
\DeclareFieldFormat[article,misc,thesis,incollection,online]{title}{\usebibmacro{string+doiurl}{\mkbibquote{#1}}}

\DeclareFieldFormat{arxiv}{%
	arXiv\addcolon\space 
	\ifhyperref
	{\href{https://arxiv.org/abs/#1}{\nolinkurl{#1}}}
	{\nolinkurl{#1}}%
}
\DeclareFieldFormat{eprint:arxiv}{
	arXiv\addcolon\space
	\ifhyperref
	{\href{https://arxiv.org/abs/#1}{\nolinkurl{#1}}}
	{\nolinkurl{#1}}%
}
\DeclareFieldFormat{rims}{%
	\ifhyperref
	{\href{http://www.kurims.kyoto-u.ac.jp/preprint/file/RIMS#1.pdf}{RIMS-\nolinkurl{#1}}}
	{\nolinkurl{#1}}%
}
\DeclareFieldFormat{mrnumber}{%
	\ifhyperref
	{\href{http://www.ams.org/mathscinet-getitem?mr=MR#1}{MR#1}}
	{MR#1}%
	}
	\DeclareFieldFormat{eprint:mrnumber}{%
\ifhyperref
{\href{http://www.ams.org/mathscinet-getitem?mr=MR#1}{MR#1}}
{MR#1}%
\addspace
eprint}
\DeclareFieldFormat{eprint}{%
\autocap{E}print\printfield{author}
\ifhyperref
{\href{#1}{available on-line}}
{#1}%
}

\renewbibmacro*{eprint}{%
\printfield{arxiv}\printfield{primaryClass}%
\newunit\newblock
\printfield{rims}%
\newunit\newblock 	 	
\printfield{mrnumber}%
\newunit\newblock
\iffieldundef{eprinttype}
{\printfield{eprint}}
}

\usepackage[pdftex,hidelinks,
hyperindex=true, bookmarks=true, bookmarksopen=true, bookmarksnumbered=true,
pdfauthor={Benjamin Collas, Murotani Takahiro, Yamaguchi Naganori},
pdftitle={Symmetries of spaces and numbers -- anabelian geometry},
pdfsubject={},
pdfstartview=FitH,colorlinks=true,
linkcolor=reference, citecolor=citation, urlcolor=e-mail
]{hyperref}

	\usepackage[noabbrev,capitalise,nameinlink]{cleveref} 

\begin{document}
\pdfbookmark{Symmetries of spaces and numbers -- anabelian geometry}{title}

\begin{abstract}	
		\emph{``Can number and geometric spaces be reconstructed from their symmetries?''}
	This question, which is at the heart of anabelian geometry, a theory built on the collaborative efforts of an international community in many variants and with the Japanese arithmetic school as a core, illustrates, in the case of a positive answer, the universality of the homotopic method in arithmetic geometry.
	
	Starting with elementary examples, we first introduce the motivations and guiding principles of the theory, then presents its most structuring results and its contemporary trends.
	
	As a result, the reader is presented with a rich and diverse landscape of mathematics, which thrives on theoretical and explicit methods, and runs from number theory to topology.
\end{abstract}

\maketitle

\setcounter{tocdepth}{3}
\tableofcontents

\newpage

\section{Reconstructions from symmetries}\label{SoNF}
\noindent At its most elementary level, for numbers and for spaces, anabelian geometry deals with properties of polynomials. We investigate how the geometric notion of symmetries can be applied to the case of numbers. The shadow of a unifying context begins to appear.

\subsection{From roots to symmetries - Galois theory}\label{Grtos}
While it is known that the field of complex numbers $\mathbb{C}$ contains the roots of all polynomials, finding explicitly\footnote{Via an algorithm that uses the four elementary operations and the extraction of roots.} such roots for a given polynomial is a more delicate task. As first noted by Évariste \textsc{Galois} ($\sim$1830), replacing roots by their symmetries -- i.e., permutations that respect the original polynomial relations -- provides deep insight on the structure of the roots.

\subsubsection{Symmetries of roots, a simple example} 
Let us first consider the polynomial $P=X^4-5X^2+6$, which factorizes into $P=(X^2-3)(X^2-2)$, and thus admits exactly the four irrational numbers $\pm\sqrt{2}$ and $\pm\sqrt{3}$ as roots. We shall describe their symmetries as permutation maps\footnote{More precisely, such permutation maps $\phi$ is defined \emph{on the field} $K=\mathbb{Q}(\sqrt{2}, \sqrt{3})$ such that (i)~for any $a\in \Q$ it holds that $\phi(a)=a$, for any $\alpha$, $\beta\in K$ it holds (ii)~that $\phi(\alpha+\beta)=\phi(\alpha)+\phi(\beta)$ and (ii')~that $\phi(\alpha\beta)=\phi(\alpha)\phi(\beta)$.} on $\{\pm\sqrt{2}, \pm\sqrt{3}\}$. A direct computation shows that there are only four possible such maps as shown in \cref{Fig:GalEx}.

\smallskip

\begin{wrapfigure}[7]{r}{.4\textwidth}\vspace{-.7em}
		\centering
			\def\arraystretch{.7}
			$\begin{array}{rcc}
				e: &\sqrt{2}\mapsto \sqrt{2},& \sqrt{3}\mapsto \sqrt{3};\\
				\varphi_1:& \sqrt{2}\mapsto-\sqrt{2},& \sqrt{3}\mapsto\sqrt{3};\\
				\varphi_2:& \sqrt{2}\mapsto\sqrt{2},& \sqrt{3}\mapsto-\sqrt{3}; \\
				\varphi_3:& \sqrt{2}\mapsto-\sqrt{2},& \sqrt{3}\mapsto-\sqrt{3}.\\
			\end{array}
			$
			\captionof{figure}{Symmetries of the field $\mathbb{Q}(\sqrt{2}, \sqrt{3})$}\label{Fig:GalEx}	
\end{wrapfigure}
Let us write $\mathbb{Q}(\sqrt{2}, \sqrt{3})$ for the field generated by $\mathbb{Q}$ and the elements of $\{\pm\sqrt{2}, \pm\sqrt{3}\}$.	The four above permutations are called the \emph{symmetries of the roots of $P$} or the \emph{symmetries of the field $\mathbb{Q}(\sqrt{2}, \sqrt{3})$}.

\medskip

For any field $K$ of this type -- called \textit{number field}, that is generated as a $\Q$-vector space by a finite number of algebraic numbers -- one obtains similarly the so called \emph{Galois group of symmetries} $\Gal(K/\mathbb{Q})$, that arises as permutations of roots of equations for the associated polynomial. In the example above, the group $\Gal(\mathbb{Q}(\sqrt{2}, \sqrt{3})/\mathbb{Q})$ is the set $\{e, \varphi_1, \varphi_2, \varphi_3\}$ with a group structure defined by $\varphi_1\circ\varphi_1=e$, $\varphi_2\circ\varphi_1=\varphi_3$, etc, to form the so called ``Klein group'', written $V_4$ -- geometrically, this group can also be seen as the group of symmetries of a rectangle.

\subsubsection{Symmetries and lattice of subfields}
The following result shows that a certain property of the field structure of such a field $K$ can be dealt with via the symmetries of its Galois group.

\smallskip

\begin{quotation} \noindent\textbf{Galois correspondence.}
{\itshape
	There is a one-to-one correspondence between the fields that are contained in $K$, and the subgroups of $\Gal(K/\mathbb{Q})$.} 
	\end{quotation}

\smallskip
	
	The example of \cref{Fig:Gcorres}, which is taken from ``A Worked out Galois Group for the Classroom''\footnote{By L.~\textsc{Halbeisen} and N.~\textsc{Hungerbühler}, in The American Mathematical Monthly, 131 (6), p.~501–510.}, to which we refer for detailed computations and notations, illustrates the potential intricacy of this correspondence in the case where the Galois group is the alternate group $A_4$ -- also the group of isometries of a tetrahedron\footnote{The reader can draw the Galois correspondence for $\mathbb{Q}(\sqrt{2},\sqrt{3})$, which is much more simple.}.
	
	\medskip
		\begin{figure*}[!htbp]
			\setcounter{figure}{0}
			\begin{subfigure}[b]{.45\linewidth}
				\includegraphics[width=\linewidth]{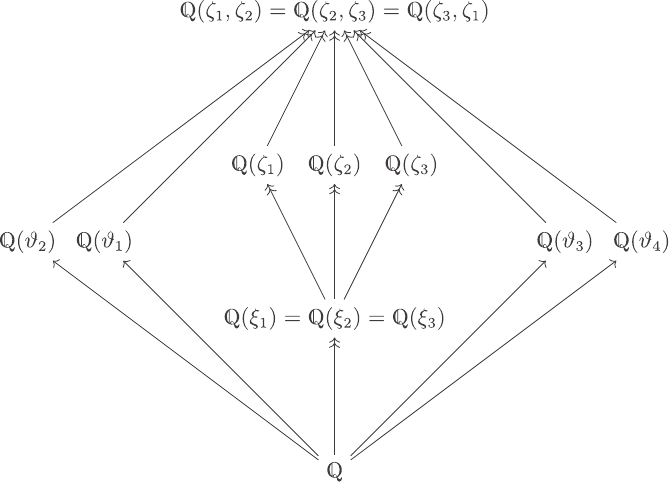}
				\caption{Fields Lattice}
			\end{subfigure} \hfill
			\begin{subfigure}[b]{.45\linewidth}
				\includegraphics[width=\linewidth]{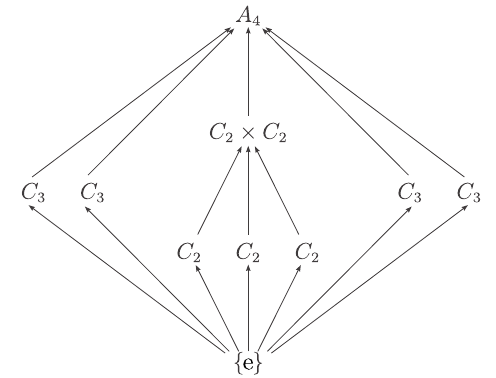}
				\caption{Groups Lattice}
			\end{subfigure}\hfill
			\captionof{figure}{Galois correspondence between the fields and groups lattices for $P=X^6-3X^2-1$}\label{Fig:Gcorres}
		\end{figure*}	
	
	By recasting an arithmetic question into a group theoretic property, this correspondence allows Galois to show that \emph{not every polynomial equation of degree greater than $4$ has computational roots.} Because the inclusion lattice is preserved, this result is also our first example\footnote{Class Field Theory -- developed by \textsc{Hilbert}, \textsc{Takagi} and \textsc{Artin} between 1880 and 1920 -- which establishes a correspondence between (the abelianization of) a Galois group and certain data attached to the number field, can be seen as another ancestor to anabelian geometry of~\cref{Sec:AG}.} of field structure information that is reconstructed from symmetries: 

\smallskip
	\begin{quotation} \noindent\textbf{First reconstruction.}
{\itshape
	The symmetries, as encoded by $\Gal(K/\mathbb{Q})$, determine the lattice of all the subfields of $K$}. 
	\end{quotation}
\smallskip
	
	We refer, once again, to the example of \cref{Fig:Gcorres}.
	
	\bigskip
	
	\noindent {\bfseries A modern glimpse.} {\itshape 
The reciprocal question, known as the \emph{Galois Inverse Problem} which was proposed by \textsc{Hilbert} in $1892$, to know if \emph{every finite group can be obtained as a Galois group over the rational numbers} is still unresolved and continues, also with its variant the ``Noether problem'', to stimulate and to shape contemporary research -- see for example Olivier \textsc{Wittenberg}'s \cite{WIT18}.
}

\subsection{...to reconstructions for numbers and spaces}\label{sec:Rfields}
The partially successful considerations of \emph{finite} symmetries for numbers and the analogy with \emph{geometric} symmetries motivate to investigate the problem of \emph{the existence of a more canonical and unifying context for reconstruction}.

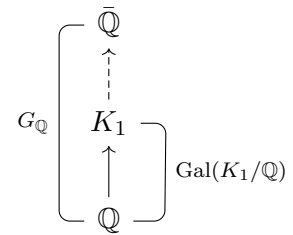
\begin{wrapfigure}[7]{r}{.2\textwidth}\vspace{-3.5em}
	\begin{tikzcd}[column sep=0pt]
		\bar{\mathbb{Q}}\arrow[no head,rounded corners, to path={--([xshift=-2ex]\tikztostart.west) -- node[left]{\scriptsize$G_\Q$}([xshift=-2ex]\tikztotarget.west)-- (\tikztotarget.west)}]{dd}\\				K_1 \arrow[no head,rounded corners, to path={--([xshift=2ex]\tikztostart.east) -- node[right]{\scriptsize$ \Gal(K_1/\mathbb{Q})$}([xshift=2.5ex]\tikztotarget.east)-- (\tikztotarget.east)}]{d}\ar[u,dashed]\\
		\mathbb{Q}\ar[u]
	\end{tikzcd}
	\caption{Tower of fields and $G_\Q$}\label{Fig:GalAbs}
\end{wrapfigure}
\subsubsection{The field structure of numbers}
In the case of numbers, one can define the field $\overline{\Q}$ of all algebraic numbers contained in $\mathbb{C}$, which thus contains all the roots of all polynomials. It contains all the towers of all finite extensions $\mathbb{Q}\subset K_1\subset \dots\subset K_n\subset\dots\subset \bar{\Q}$, where, for example, $K_1=\mathbb{Q}(\sqrt{2}, \sqrt{3})$ as in the previous section.

On the group side, one can then form, as the inverse limit of the $G(K_i/\Q)$, \emph{the absolute Galois group} of rational numbers $G_\Q=\Gal(\overline{\Q}/\Q)$, which is the seed of number theory, and whose structure\footnote{The absolute Galois group of rational numbers is an infinite topological group -- more precisely, a profinite group, i.e., obtained as inverse limit of finite groups.} is quite rich. While, in terms of elements, much remains unknown about $\Gal(\overline{\Q}/\Q)$, in terms of symmetries, an analogue of the fundamental theorem of Galois theory holds: \emph{any field $K$ between $\Q$ and $\overline{\Q}$ corresponds to one and only one (closed) subgroup of $\Gal(\overline{\Q}/\Q)$}.

\medskip

It is remarkable that, in this context, and as established in successive steps by J\"{u}rgen \textsc{Neukirch}, \textsc{Ikeda} Masatoshi, \textsc{Komatsu} Keiichi, \textsc{Iwasawa} Kenkichi, and \textsc{Uchida} K\^{o}ji ($\sim$1970), one can go beyond the reconstruction of the field lattice of the previous section:

\smallskip

\begin{quotation} \noindent\textbf{Galois reconstruction of number fields.}
\itshape
There exists an \emph{algorithm} which, starting from a group $G$ of the type of an absolute Galois group of a number field, gives the reconstruction of the number field~$K$.
\end{quotation}

\smallskip

Here ``algorithm'' must be understood as 
\smallskip

\begin{quotation}
\noindent\emph{``a definite and functorial process based uniquely on the topological group structure of the input $G$''}.
\end{quotation}

\smallskip

For example, in the case where $K$ is a $p$-adic field, as in \cref{{ft:pAdic}}, the Artin reciprocity map of algebraic number theory provides the isomorphism
\[
G^\mathrm{ab}\simeq \widehat{\mathbb{Z}}\times \mathbb{Z}/(\mathbf{p_K^{f_K}}-1)\mathbb{Z} \times \mathbb{Z}/\mathbf{p_K}^b\mathbb{Z}\times \mathbb{Z}_{p_K}^\mathbf{d_K},
\]
which from $G$ recovers the ramification index $e_K=d_K/f_K$, then the inertia $I_K\simeq\cap N\triangleleft G$ (where the intersection is taken over the $e_K$-cofinite normal subgroups), to recover \emph{the multiplicative monoid} $(K^\times,\boxtimes)$ from $K^\times \simeq G^\mathrm{ab}\times_{G/I} \mathbb{Z}$. The algebraic closure $(\bar{K}^\times,\boxtimes)$ follows by the functoriality of local class field under ``Verlagerung'' maps and direct limit. A similar process reconstructs the \emph{additive topological module\footnote{While \emph{the final topology} induced by the intermediate field extensions is group-theoretic, the usual $p$-adic topology on $K$ is not -- see \cite{HOS22} Rem.~4.3.2.}} $(\bar{K},\boxplus)$ -- see \cite{HOS22}~\S~3-4.

\medskip

\begin{wrapfigure}[15]{r}{.5\textwidth}\vspace{-1em}
	\centering
	\includegraphics[width=.8\linewidth]{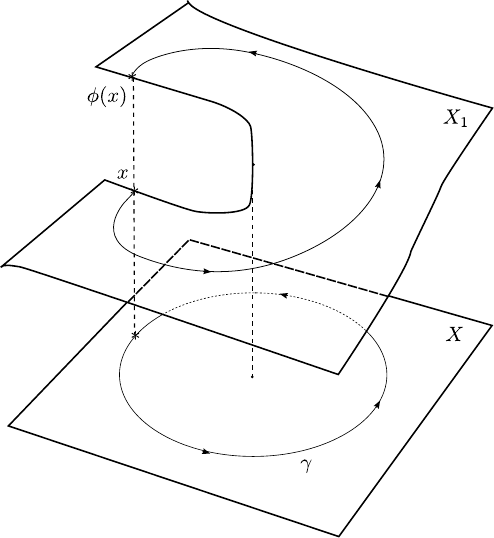}
	\caption{Symmetries: a loop and a covering.}\label{Fig:PathCover}
\end{wrapfigure}
The Galois symmetries thus encode all the information of $K$ and its field structure. Note that this algorithmic result, that involves only one group (and not the comparison of two), is a contemporary refinement of \textsc{Hoshi} Yuichiro in \cite{HOS22} of the original result.

\subsubsection{A geometric recasting}\label{sec:HTMstow}
Let us see how the group of Galois symmetries for numbers admits an equivalent for geometric spaces\footnote{Or \emph{manifolds}, spaces, in our case over $\mathbb{C}$, which can be of any dimension and are obtained as gluing pieces of (our) usual euclidean spaces, and which can have holes and doughnut-like shapes.}, that is the group of loops on the space.

\medskip

Loops on a topological space $X$ (over $\mathbb{C}$) are continuously deformable closed curves (with same starting and ending point $*$), which can be composed by concatenation to form a group called the \emph{topological fundamental group}, written $\pi_1^\mathrm{top}(X,*)$. If we consider a (finite) covering of $X$ by another manifold $X_1$, as in \cref{Fig:PathCover}, we notice that a loop $\gamma$ on $X$ defines, by sending the starting point $x$ to the end point $\phi(x)$, a transformation $\phi$ of $X_1$. The corresponding group of transformations $\Aut(X_1/X)$ is the analog of Galois symmetries $\Gal(K_1/\mathbb{Q})$ for spaces.

\smallskip

Pushing the Galois analogy further, this raises the question to know \emph{if, in turn, geometric spaces can also be determined by their symmetries}. In the case of manifold that are locally saddle-like, or \emph{hyperbolic}, one has:

\smallskip

\begin{quotation} \noindent\textbf{Mostow's rigidity theorem~($\sim$1968).}
\itshape
The topological fundamental group of a \textit{hyperbolic manifold} of dimension greater than $2$ completely determines the manifold.
\end{quotation}

\smallskip

At this stage, we thus obtain two separated reconstructions from symmetries: \emph{in arithmetic} with the Neukirch-Ikeda-Iwasawa-Uchida+ theorem , and \emph{in geometry} with Mostow's rigidity theorem. As we will see in the next section, the unifying context for our reconstruction is provided by Grothendieck's arithmetic geometry.

\bigskip

\noindent {\bfseries A modern glimpse.} {\itshape  An extension of the original (non-algorithmic)  Neukirch-Ikeda-Iwasawa-Uchida result was extended by Florian \textsc{Pop} (1994). More recently, it was shown by \textsc{Tamagawa} Akio and Mohamed \textsc{Saïdi} that only a smaller portion of the Galois group is sufficient to reconstruct the field -- the ``$m$-step solvable Neukirch-Uchida theorem'' -- see \cite{SAITAM22} and also \cite{POP21}.
}

\section{The universality of homotopic arithmetic geometry}\label{Sec:AG}
\noindent It follows the first volume of the ``Séminaire de Géométrie Algébrique'' (IHES, 1960) that the unifying context of arithmetic and geometry, or ``homotopic arithmetic geometry'', is provided by algebraic varieties and their étale fundamental groups. 

\subsection{Étale reconstructions for algebraic varieties}
\begin{figure}[ht]\centering
	\includegraphics[width=.8\linewidth]{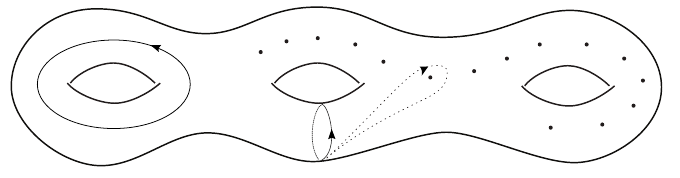}
	\caption{A hyperbolic curve of genus $3$ with $n$ marked points and a (dotted) cuspidal generator}\label{fig:RS}
\end{figure}

An algebraic variety $X$ defined over a field $K$ is the analog of a manifold obtained by patching the zero loci of polynomial equations with coefficients in $K$ -- we refer to \cref{fig:RS} for an example of a smooth curve\footnote{A one-dimensional complex variety, or curve, appears as a two-dimensional real variety, or surface.} over $K=\mathbb{C}$. 

\medskip
\begin{wrapfigure}[10]{r}{.35\textwidth}\vspace{-4em}
					\centering	
					\includegraphics[width=\linewidth]{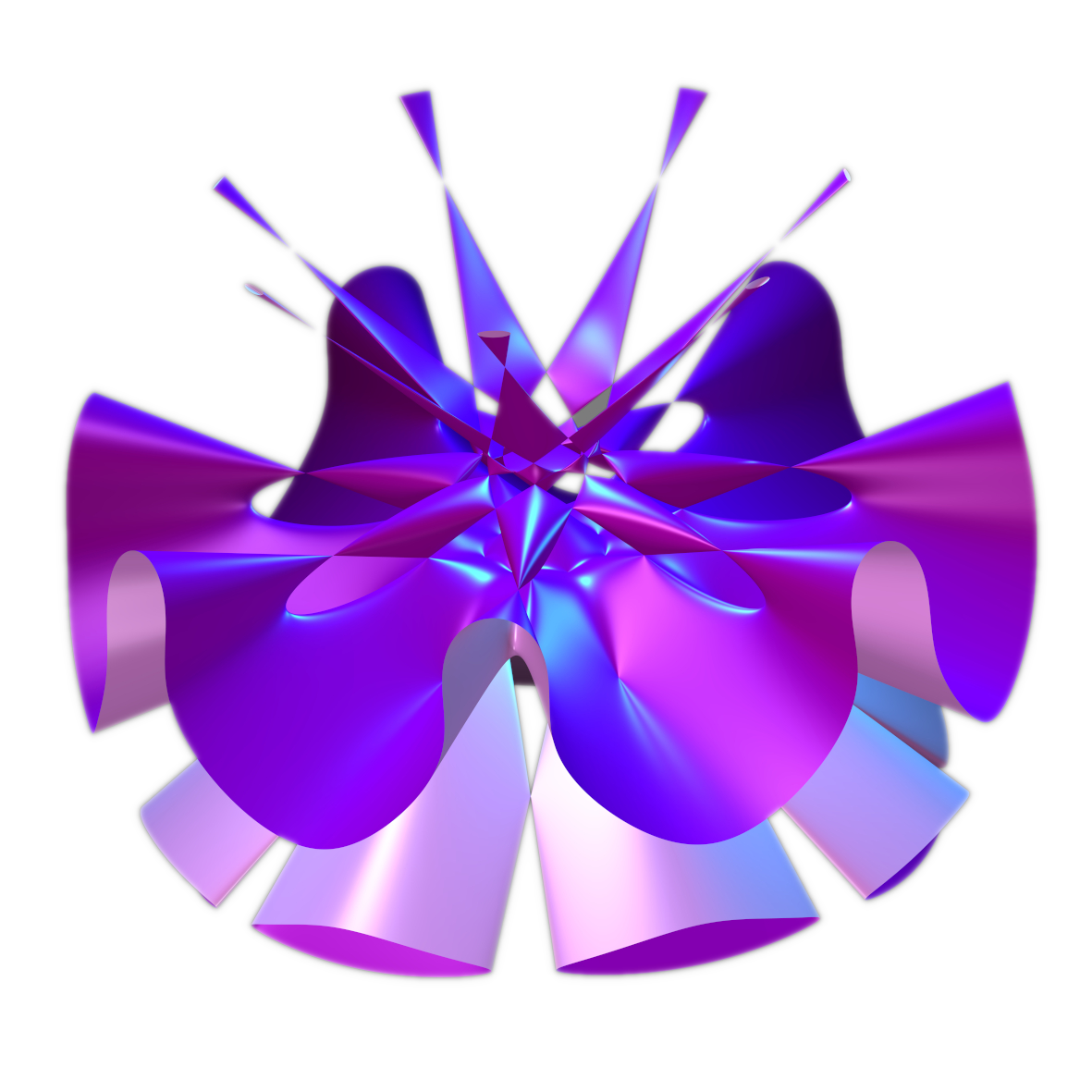}
					\caption{An algebraic surface of degree $7$}\label{Fig:AlgSurf}
		\end{wrapfigure}
Similarly to the previous section, the \emph{étale fundamental group} $\pi_1^\mathrm{et}(X,*)$ encodes the transformations -- or étale symmetries -- of certain types of coverings of $X$. The étale fundamental group is a profinite topological group, which, since one has the following identifications 
\begin{itemize}
	\item for a manifold $X$ over $\mathbb{C}$: $\pi_1^\mathrm{et}(X,*)=\hat{\pi}_1^\mathrm{top}(X,*)$,
	\item for a number field: $\pi_1^\mathrm{et}(\mathrm{Spec}\,K)=\Gal(\bar K/K)$,
\end{itemize}
generalizes both the absolute Galois group and the topological fundamental group.

\medskip

\subsubsection{When the arithmetic meets the geometry}\label{sec:ArMeetGeo}
Let $K$ be a number field and $X$ be an algebraic variety over $K$. By standing in the following exact sequence -- where $X^\mathrm{an}$ denotes the complex variety defined by the polynomial equations of $X$ with solutions taken in $\mathbb{C}$ instead of $K$: 
\begin{equation*}\label{eq:FES}\tag{FES}
1\to \underbrace{\hat{\pi}_{1}^\mathrm{top}(X^{\mathrm{an}},*)}_\text{Geometry}\to \pi_{1}^\mathrm{et}(X,*)\to \underbrace{\mathrm{Gal}(\bar K/K)}_\text{Number theory}\to 1,
\end{equation*}
the étale fundamental group \emph{intermingles number theory and geometry}.

\medskip
The sequence (FES) above further defines some actions of the absolute Galois group on (a completed version of) the topological fundamental group. 

\medskip

Let us recall the classical example of $X=\mathbb{P}^1\setminus \{0,\infty\}$ over $K=\Q$, where, topologically, $X^\mathrm{an}$ is the complex plane with the point $0$ removed, so that the topological generating loop $x$ around $0$ of  is a generator of the profinite $\widehat{\pi}_{1}^\mathrm{top}(X^{\mathrm{an}},*)$. For each finite $N$-cover of the latter, the loop $x$ identifies with a $N$-root of unity $\zeta_N$, on which $\sigma\in\Gq$ acts by a certain $N$th-power $\zeta_N\mapsto \zeta_N^{\chi_N(\sigma)}$. By taking the inverse limit $\widehat{\pi}_{1}^\mathrm{top}(X^{\mathrm{an}},*)\simeq \varprojlim_N \mathbb{Z}/N\mathbb{Z}=\widehat{\mathbb{Z}}$, one thus obtain the $\Gq$-action given by the cyclotomic character $\chi=\varprojlim_N \chi_N$ (also an element of $\widehat{\mathbb{Z}}$):
\begin{equation}\tag{cyc}\label{Eq:cycAct}
\sigma.x=x^{\chi(\sigma)}\in \pi_1^\mathrm{et}(\mathbb{P}^1_{\bar{\Q}}\setminus \{0,\infty\},*)
\end{equation}	

\begin{wrapfigure}[13]{r}{.5\textwidth}\vspace{-1em}
	\centering
		\begin{tikzpicture}[every braid/.style={
				ultra thick,
				braid/strand 1/.style=Mycolor1,
				braid/strand 2/.style=Mycolor2,
				braid/strand 3/.style=Mycolor4,
				braid/anchor=center,
			}]
			\pic[
			braid/.cd,width=1cm,
			] {braid={s_1 s_1 s_2 s_3}};
		\end{tikzpicture}
		\caption{Under $\Gq$, this braid is sent to a conjugate and $\chi$-power of itself as in \eqref{Eq:cycAct}.}\label{Fig:Brd}	
\end{wrapfigure}

The action of $\Gq$ thus appears to be computable, also on similar spaces, such as the moduli spaces of curves $\mathcal{M}_{g,[n]}$ or the configuration spaces of $n$ points on (hyperbolic) curves $\mathrm{Conf}_n(X)$: the local action above is transported along analytic paths from ``$0$'' to every missing point or component; similar formulas to \eqref{Eq:cycAct} ensue, for example on Artin braid generators and certain braid words such as in \cref{Fig:Brd}.

\medskip

This constitutes the first stage of encoding $\Gq$ via the geometric combinatoric of spaces -- the so called Galois-Teichmüller theory\footnote{Which provides theoretical and combinatorial inputs for the reconstruction of varieties from their étale symmetries.}, see \cite{OES03} for an introduction.			

\subsubsection{Anabelian arithmetic geometry, then...}
Anabelian\footnote{Contrary to a first idea,  ``Anabelian'' does not stand for ``ana-belian'', from the Greek ``ana=anew'', but for ``an-abelian'' with ``an=without'' and ``abelian=commutative'' -- from the mathematician Abel. Anabelian geometry deals with objects whose fundamental group lacks commutativity.} arithmetic geometry deals with the inverse process of the previous section, that is the question of the reconstruction of spaces from exact sequences of the type of (FES). This kind of geometry emanates from a 1987 insight of Alexander \textsc{Grothendieck}\footnote{As envisioned in his foundational letter to Gerd \textsc{Faltings} in 1983, which, as it since appeared and is emphasized in \cite{BOYRT25}, is ``very roughly in the right direction, but not really in the right direction''.}, which, in its original form, states that

\smallskip

\begin{quotation}\itshape
\noindent Any isomorphism between the étale fundamental groups of two hyperbolic curves comes from a unique isomorphism between the varieties themselves.
\end{quotation} 

\smallskip

This conjecture has been resolved by successive progress, which each introduces decisive and lasting arithmetic insights:

\medskip
\begin{itemize}
	\item[(R1)]~\textbf{Genus $0$ curves over number fields}, by \textsc{Nakamura} Hiroaki (1991).
A \emph{group-theoretic} property of $\pi_1^\mathrm{et}(X)$ leads to a linear Galois action (i.e., on finite dimensional vector spaces), for which Deligne theory of weights applied to the cusps-abelianization sequence 
\begin{equation*}
...\to \overbrace{\underbrace{\bigoplus\widehat{\mathbb{Z}}}_\text{weight $-2$}}^\text{cusps}\to \pi_{1}^\mathrm{et}(X,*)^\mathrm{ab}\to \overbrace{\underbrace{\pi_{1}^\mathrm{et}(\bar{X},*)^\mathrm{ab}}_\text{weight $-1$}}^\text{abelian}\to 1,
\end{equation*}
where $\bar{X}$ denotes a compactification, group-theoretically characterizes the genus and number $n$ of marked points, or cusps, of $X$. One further obtains \emph{for any genus} the reconstruction of the inertia groups $I_x$ (or loops in $X$ around the marked points) which are generators of the étale fundamental group, see \cite{NaK90}.

\smallskip

\item[(R2)]~\textbf{Affine curves, any genus, and over number fields}, by \textsc{Tamagawa} Akio~(1997), see \cite{Tam97}. By group theoretically counting rational points ``over $X$'', one obtains some limit points $x_\infty\in X$ to which (R1) can be applied to reconstruct the decomposition groups of closed points . The reconstruction of the field $K(X)$ of functions over the space\footnote{\label{FN:AG}In Grothendieck's geometry, algebraic spaces are not only defined by their points but also by a certain \emph{ring} of functions $(\mathcal O_X,\boxtimes,\boxplus)$, which allow introducing various topologies in relation with the arithmetico-geometric problem at play.} $X$ is dealt with in terms of Kummer theory and étale cohomology \emph{groups}; one reconstructs \emph{the multiplicative structure first}, then the additive one.

This approach leads to the general notion of ``Kummer-faithfulness'' for the base field $K$.

\smallskip

\item[(R3)]~\textbf{Any curves and genus, and over $p$-adic fields\footnote{\label{ft:pAdic}A $p$-adic field $K/\Q_p$, that generalizes number fields (in the sense that any number field can be embedded in such a field, see Cassel's $p$-adic embedding Theorem), can be thought of as a neighborhood of a prime $p$ in $\Q$. It must be noted that Grothendieck's original conjecture for $K$ a $p$-adic field does not hold (but a refined version does): there exists non-isomorphic $p$-adic fields with isomorphic Galois groups.}}, by \textsc{Mochizuki} Shinichi~(1999).
The étale coverings of $\pi_1^\mathrm{et}(X)$ are dealts with via line bundle techniques, while $p$-adic Hodge theory provides a container for $X$ made of differential forms: $\mathbb{P}(\Gamma_X)\hookleftarrow X$, see \cite{MOC99}.
\end{itemize}


At this stage, anabelian geometry \emph{exploits or reformulates some classical techniques from}
\begin{itemize}
	\item[(a)]~\emph{Number theory:} local and global class field theory, for (R2); and 
	\item[(b)]~\emph{Algebraic geometry: Deligne theory of weights for (R1), an anabelian good reduction criterion à la Néron–Ogg–Shafarevich, and the Lefschetz trace formula to isolate rational points in covers as in (R2)}. 
\end{itemize}	

Let us refer to the still pertinent and illuminating survey \cite{NTM01} for more details and additional comments\footnote{E.g., the ``Galois'' or birational variant of Grothendieck's conjecture, that deals with spaces of the form $X=K(x_1,\dots,x_n)$ with $K=\Q$ or $\mathbb{F}_p$ was solved by \textsc{Pop} (1994).}, completed with \cite{HOS25}. 

It also \emph{extends fundamental previous results} such as: the Neukirch-Ikeda-Iwasawa-Uchida theorem (from dimension $0$, a field, to dimension $1$, a curve) or Faltings' isogeny theorem (from abelian to non-abelian groups) which is a key step in establishing the finiteness of rational points on curves (the Mordell conjecture) -- see~\cite{NTM01}  \S~1.2 and \S~4.1.

\medskip

Anabelian geometry has furthermore crystallized into ``group-theoretic reconstructions'' and into the \emph{contemporary mono-anabelian reconstruction} (see also \cite{HOS22} \S~Introduction), a consideration which departs from Grothendieck's original insight by involving only one geometric space.

\subsubsection{Towards higher ``dimensions''?} 
Before dwelling into more structuring contemporary considerations, the attentive reader naturally, and more prosaically, has certainly wondered: \emph{``After the dimension $0$ and the dimension $1$... is there, similarly to Mostow's theorem, some anabelian reconstructions in higher dimension\footnote{Two other directions are given by the investigation of arithmetic phenomenons (a)~in higher homotopy, i.e., when taking the $\pi_N^\mathrm{et}(X)$, $N\geq 1$, into account, and (b)~in higher categories, i.e., when exploiting the stack structure of the spaces -- we refer respectively to Schmidt and Stix' \cite{SS18} and to the first author's \cite{COL21}.}?''} 

\medskip

One difficulty comes here from the so-called ``purity property'' of the étale fundamental group, which does not see any information in subspaces of small dimension (i.e, of codimension at least $2$) -- low-dimensional group reconstructions cannot be glued into higher dimensional reconstructions.

\begin{figure}[!h]\centering
$X=X_n\to X_{n-1}\to\dots\to X_1\to X_0=\mathrm{Spec}\ K$
\caption{A polycurve $X$ over the field $K$ (where $X_{i+1}\to X_i$ are families of hyperbolic curves)}\label{Fig:PolyC}
\end{figure}

Building on the previous techniques of (R1)-(R3), and because \emph{anabelian reconstructions spread well in families of curves} as in Figure~\ref{Fig:PolyC}, one can nevertheless establish, in the case of $g=0$, the \emph{anabelianity of $\mathrm{Conf}_n(C)$ and $\mathcal{M}_{0,n}$} --  put $X_i=\mathrm{Conf}_i(C)$ in the above -- see \textsc{Nakamura} Hiroaki, \textsc{Takao} Naotake, and \textsc{Tamagawa} Akio Theorems~C of~\cite{IN97} -- then \cite{AbsTopII} Cor.~1.11 for any genus. Lie algebra techniques recently provided an \emph{explicit algorithm to group theoretically reconstruct} the genus, the number of points and the dimension of $\mathrm{Conf}_n(C)$ -- by \textsc{Sawada} Koichiro in~\cite{Saw18}.

\medskip

This technique of polycurves (latter of quasi-tripods, see below) was even more recently used to establish that \emph{(smooth) algebraic varieties of any dimension possess a fundamental system of neighborhoods that are anabelian}, another conjecture of {Grothendieck} -- see \textsc{Hoshi}'s \cite{HOS21}, also \textsc{Schmidt}-\textsc{Stix}'s \cite{SS18} with different techniques.

\medskip

\noindent\textbf{A Galois-homotopy remark.~}{\itshape	
A version of the Inverse Galois Problem purely for $\pi_1^\mathrm{et}(X^{an})$ does exist; the so called \emph{Abhyankar’s conjecture} was finally solved for curves by \textsc{Harbater} in 1997.
Michael \textsc{Frieds}'s geometrization of the inverse Galois problem in terms of $G$-covers of genus $0$ curves -- then followed by Pierre \textsc{Dèbes}, Moshe \textsc{Jarden}, and others -- brought the IGP closer to anabelian arithmetic geometry. An unifying context is finally given by considering $G$ as a \emph{hidden stack symmetry} group of curves, see \cite{COL21}.
}

\subsection{... and now: contemporary anabelian developments}
Beyond Grothendieck's (relative) initial vision, anabelian geometry nowadays focuses on \emph{establishing and exploiting the universality of arithmetic homotopy}\footnote{As opposed, for example, to number theory (that is too rigid and ``does not see enough'') or complex algebraic geometry (that is too flexible and ``sees too much'').}, that is, how spaces and numbers can be studied via the canonical group-theoretic properties of their étale fundamental groups. A panorama of recent progress suggests the two following questions:

\smallskip

\begin{quotation}
\noindent{\bfseries Universality of anabelian arithmetic geometry.}
\itshape
\begin{enumerate}
\item \label{it:1} Which kind of new insight, in number theory and algebraic geometry, is given by the anabelian homotopic method?
\item \label{it:2} Does there exist a group-theoretic algorithm, which from a group that is of étale fundamental group type, reconstructs the original space?
\end{enumerate}
\end{quotation}

\smallskip

These structuring questions illustrate, respectively, \emph{the ubiquity} and \emph{the canonicity} of the homotopic method in arithmetic geometry. In return, they both lead to beautiful mathematical insights and interactions\footnote{We refer to \cite{BOYRT25} for a lively exchange between \textsc{Hoshi}, \textsc{Mochizuki}, \textsc{Tamagawa}, and the first author on this topic.}.


\subsubsection{Canonicity: algorithms and absolute mono-reconstructions}
The canonicity of anabelian geometry is mostly expressed in the consideration of algorithms that, in the construction steps, eliminate any choice and rely on group-theoretic arguments only, and in the consideration of categorical contexts. 

\medskip

\cref{algo:anab} provides a sketch of such an algorithm\footnote{We must insist that this algorithm applies to very specific \emph{types} of curves only, and that \emph{the algorithm}, not the resulting isomorphism, is important.} that reconstructs the function field $K(X)$ and the base field $K$ of a hyperbolic curve of strict Bely\v{\i} type (see \cref{Fig:StcB} for a definition, and Theorem~1.11 in \cite{AbsTopIII}).

\begin{algorithm}\small
\caption{Mono-anabelian reconstruction for a hyperbolic curve of strict Bely\v{\i} type}\label{algo:anab}
\begin{algorithmic}[1]	
\State \textbf{Geometric group:} As the maximal finitely generated closed normal subgroup of $\Pi\rightsquigarrow\Delta$
\State \textbf{Genus and points:} As the dimension of vector spaces in weighted étale cohomology $\Delta\rightsquigarrow (g,n)$
\State \textbf{Inertia groups:} By Bely\v{\i} cuspidalization, the decomposition, then inertia groups $(\Delta, n)\rightsquigarrow D_x\text{ and } I_x$
\State \textbf{Multiplicative  monoid:} By Kummer theory $\rightsquigarrow$ multiplicative monoid $(K(X)^\times,\boxtimes)$
\State \textbf{The (function field of) $X$:} By Uchida Lemma and divisors properties $\rightsquigarrow (K(X),\boxtimes,\boxplus)$
\State \textbf{The base field of $X$:} $K(X)\rightsquigarrow K_X$
\end{algorithmic}	
\hrule\vspace{.2em}
{\footnotesize Here $\Pi$ (resp. $\Delta$) denotes a group isomorphic to a certain $\pi_1^\mathrm{et}(X,*)$  (resp. $\hat{\pi}_1^\mathrm{top}(X^\mathrm{an},*)$).}
\end{algorithm}			 	

As always, the attentive reader will have noticed that \cref{algo:anab} doesn't reconstruct the curve $X\to\mathrm{Spec}\ K$ from the portion $\pi_{1}^\mathrm{et}(X,*)\to \mathrm{Gal}(\bar K/K)$ of the (FES), but rather reconstructs the curve $X$ \emph{above} $K$ from $\Pi=\pi_{1}^\mathrm{et}(X,*)$ only. This \emph{absolute} reconstruction (see \cite{MOC12}) allows the structure morphism $X\to \mathrm{Spec}\ K$ to vary -- a process tantamount to allow the variation of the coefficients in the polynomials defining $X$, and thus a more canonical definition of the geometric locus associated to $X$.

\subsubsection{Canonicity: new structures appears}
This approach in particular reveals the importance of two types of structures.

\medskip

\noindent (S1)~\textbf{Rings become multiplicative and additive monoids}. Grothendieck's algebraic geometry, that is the modelization of geometric spaces $X$ by polynomials, relies on the ring structure of the function \emph{ring} $(\mathcal O_X,\boxtimes,\boxplus)$ of $X$ (see \cref{FN:AG}). 

On the other hand, in the reconstruction of a $p$-adic field $K$ (that is in dimension $0$ already), appears some objects related to $K$, whose ring structure cannot be anabelianly reconstructed. That is the case of the field $K$ itself as in \cref{ft:pAdic}, and also, for example, of the units of the ring of integers $(\mathcal{O}^\times_{\bar K},\boxtimes,\boxplus)$, whose $\boxtimes$-multiplicative and $\boxplus$-additive monoid structures can nevertheless separately be recovered -- see~\cite{HOS22}. 

\medskip

\begin{wrapfigure}[7]{r}{.4\textwidth}\centering\vspace{-1.2em}
\small
\begin{tikzcd}[row sep=4pt, column sep=15pt,/tikz/row 2/.append style={row sep=20pt}, ampersand replacement=\&]
V \ar[ddr, twoheadrightarrow, bend right, "f.et."'] \ar[dr,twoheadrightarrow, dotted,  "fin.","et."']\& \& U_X\ar[ddl,hook, dotted, bend left,"op."] \ar[dl,twoheadrightarrow, dotted,  "\exists\ fin."', "et."]\\
\&\mathbb{P}^1\setminus D  \& \\
\& X \&	\\[-7pt]
\end{tikzcd}
\caption{Curve $X$ of strict Bely\v{\i} type}\label{Fig:StcB}
\end{wrapfigure}
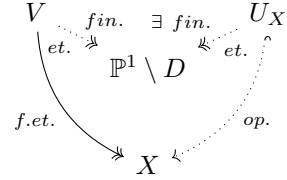
This situation, that also exists for a certain type of genus $1$ hyperbolic curves, results in a (weaker) \emph{geometry of monoids}, where \textbf{(a)}~$\boxtimes$ and $\boxplus$-monoidal structures can vary independently, and where \textbf{(b)}~``anabelian reconstructible'' objects can be used as shared containers to compare the variations -- see ibid.

%
%
%

\medskip

\noindent (S2)~\textbf{New categories of curves: curves of strict Bely\v{\i} type and quasi-tripods}.	
{\itshape How to ensure the \emph{mono}-anabelian reconstruction of the inertia group of (R1)~above?} The answer is given by the \emph{categorical context} of curves of strict Bely\v{\i} type -- see Corollary 3.7 of \cite{AbsTopII} for a complete algorithm.

For such a curve $X$ to stand in a diagram as in \cref{Fig:StcB} results in a finite sequence of étale quotients and open immersions
\[
X\rightsquigarrow V \overset{1...n}{\rightsquigarrow} \mathbb{P}^1\setminus D\overset{1...r+m}{\rightsquigarrow} U_X\rightsquigarrow X
\]
which says that ``cusps'' (as on the right side of \cref{Fig:StcB}) can be controlled in terms of finite étale covers (as on the left side of \cref{Fig:StcB}).

\medskip

The definition of a quasi-tripod curve $X$, which leads to the existence of fundamental system of absolute anabelian neighborhoods for any algebraic variety (and the anabelianity of $\mathrm{Conf}_n(C)$ in any genus, profinite case), is based on similar considerations:
\[
X\rightsquigarrow X_1 \leftrightsquigarrow ... \leftrightsquigarrow X_n\leftrightsquigarrow  \mathbb{P}^1\setminus\{0,1,\infty\},
\]
where one considers finite étale covers and open immersions -- we refer to \cite{Tripod} for anabelian reconstruction results.

\subsubsection{Canonicity: categorification}
The \textbf{combinatorial anabelian geometry} of \textsc{Hoshi} and \textsc{Mochizuki}, which deserves a survey on its own -- for a first contact, we recommend \cite{SM10} (introduction included) and \cite{CAGG07}, also \cite{HOS25}~\S7-8 for an overview of recent results -- is built by \emph{attaching certain Galois categories} on each smooth components of degenerated curves. It can be seen as a functorial\footnote{Also as an arithmetic version of the pants (or ``tripods'') approach of Hatcher-Thurston to obtain a presentation of the mapping class group of surfaces.} version of the Galois-Teichmüller theory for braids groups of \cref{Sec:AG}. 

By controlling the inertia groups in fibration (as for (R1) and in Figure~\ref{Fig:PolyC}), and by group theoretically synchronizing the various projective lines (or ``tripods'') of the configuration, it stands at the source of anabelian results in higher dimensions. 

From this follows \emph{the anabelianity of Galois-Teichmüller theory}, as in \cite{HMM22}, then a promising combinatorial model $\bar{\Q}^\mathrm{BGT}$ of the algebraic closure $\bar{\Q}$ of the rational numbers  -- see \cite{HMT20} and more generally \textsc{Tsujimura} Shota et al. in~\cite{TSU23}.

\section{Ubiquity: exploration and new frontier}
Having identified new arithmetico-geometric invariants, the new role of certain structures, and new techniques, anabelian geometry is now turning around by providing new insights in arithmetic and geometry\footnote{``Arithmetic Geometry'' exists in a Diophantine variant, following Serge \textsc{Lang}, a ``(co)homologic'' variant, following Deligne's IHES seminar with a focus on linear Galois actions, and a ``homotopic'' variant, following Grothendieck's ``Récoltes et semailles'', as in this manuscript (see \S2.8 and \S~2.10 ibid).}.

\subsection{Ramifications and explicit methods}
Let us briefly mention a few among the many recent developments, and refer the reader to the corresponding author's reports in the already cited ``Oberwolfach Report'' of 2023 for details and original references. 

\subsubsection{Algebraic geometry}
\textsc{Tamagawa}'s \textbf{resolution of non-singularities} provides the existence of certain type of models for curves. Intersecting with Berkovich $p$-adic analytic geometry, it further establishes Grothendieck's \emph{absolute conjecture over $p$-adic fields} -- see~\cite{MT23} and more generally \cite{LEP23}.

For \textbf{elliptic curves}, an absolute mono-anabelian reconstruction algorithm of the theta function, as in Theorem~1.6 of \cite{EtTh}-- which classicaly provides an analytic projective model of the curve -- exploits the weaker monoid structures of (S1) above (\textsc{Mochizuki}).

\subsubsection{Algebraic number theory}
It follows \textbf{Ihara's program} that arithmetic homotopic data of (FES) can replace classical methods of algebraic number theory: \textbf{(a)}~for the study of certain Galois extensions as in the ``When does the mountain meets the Heaven?'' journey pursued by \textsc{Rasmussen}-\textsc{Tamagawa} since 2017 -- see \cite{MR24}, and as in Oda's prediction \cite{CP25}, or \textbf{(b)}~for revisiting \textbf{Greenberg's program},  as by Rachel \textsc{Pries} \cite{PRIES23}. 

One further observes new connections (e.g., in relation to fields properties, such as the anabelian ``Kummer faithfulness'' of (R2)) with \textbf{classical Galois theory \emph{à la} Field Arithmetic} -- see the second author's \cite{MUR23}, \textsc{Sawada} et al. in \cite{MST24} with a refinement of ``indecomposibility'' (a property shared by Galois groups \cite{MT22} and GT \cite{M18-GT}), and recent project of \textsc{Taguchi} Yuichiro on the Mordell-Weil group \cite{AT24}.

\subsubsection{Low-dimensional topology}
The Teichmüller mapping class groups and braid groups -- and more recently their categorical version, the little $2$-discs operads $E_2$, see \textsc{Fresse}-\textsc{Horel} as in \cite{Fresse19-cs} -- naturally appears as fundamental groups of the $\mathcal{M}_{g,[n]}$. They provide, for example, an \emph{arithmetic obstruction} to the Johnson homomorphism's surjectivity that completes Morita's -- see \textsc{Nakamura} as in \cite{Nak96}.

\medskip

Many of these approaches are explicit by involving computations on matrices, braids, or Lie algebras\footnote{For a \texttt{Lean} library for anabelian geometry, see Topaz' project (in progress).}. 

\subsection{A homotopy-homology frontier?}
(Co)Homology groups are abelian, the first étale homotopy one is not. In contradiction with Grothendieck's original vision, the \textbf{nearly-abelian or minimalistic program} -- or the $m$-step solvable Grothendieck conjecture, which implies the original one -- investigates how ``anabelianly big'' the étale fundamental group must be to reconstruct the spaces. 

As it appears, the whole fundamental group is not always necessary, \emph{a nearly-abelian quotient sometimes suffices} -- see \textsc{Saïdi-Tamagawa} for Galois groups, the geometric version of the third author \cite{YAM24}, the \textsc{Pop}'s Neukirch-Uchida \cite{POP21}, and the Galois-Teichmüller approach of Adam \textsc{Topaz}~\cite{PT21}.

We also refer to the Deligne-Ihara Lie algebra in genus $1$ of \textsc{Ishii} Shun \cite{ISHII25} with potential link with \textbf{motivic theory}, the universal cohomology theory (see also \cite{COL21}).

\medskip 

As presented\footnote{Except the reconstruction over $\bar{\mathbb{F}}_p$ initiated by \textsc{Tamagawa} in \cite{Tam90} (and more recently, the extension of \cite{Yang22}), which, following \textsc{Hoshi} and Yu \textsc{Yang}, now provides an independent proof of (R3) -- paper in preparation.} in these notes, anabelian geometry has thus grown, in techniques and results, beyond Grothendieck's original insight to indicates new research frontiers in arithmetic geometry; we refer the interested reader to the \S~Introduction of the ``Oberwolfach report'' 2023, which is the result of a renewed collaborative effort of the Arithmetic and Homotopic Galois Theory community.

\bigskip

\noindent{\itshape \S~Acknowledgments.}
{\small This manuscript is part of the France-Japan AHGT international research network supported by the Research Institute for Mathematical Sciences of Kyoto University and CNRS. Second author, resp. third author, is supported by the (JSPS) KAKENHI Grant Numbers 22J00022 and 24K16890, resp. Number 23KJ0881. The first author expresses his gratitude to Kiran S.~\textsc{Kedlaya} and Dinesh \textsc{Thakur} for comments on a previous version of this manuscript.}

\medskip

\noindent{\small\emph{Versions history.} Jan. 2025: submitted; Apr. 2025: +5 pages; July 2025: Extended bibilography.}

\medskip
	
	{\scshape\raggedleft\large References\par}
	\addcontentsline{toc}{section}{References}
	\begin{multicols}{2}
	\printbibliography[heading=none] 
	\end{multicols}

{\centering
	$\ast\quad \ast \quad \ast$
	\par}

\medskip

{\small
	\noindent Benjamin \textsc{Collas} is researcher in homotopic arithmetic geometry at the \emph{Research Institute for Mathematical Sciences of Kyoto University}, Japan \textbullet Email: bcollas@kurims.kyoto-u.ac.jp \textbullet \url{https://www.kurims.kyoto-u.ac.jp/~bcollas/}
	
	\medskip
	
	\noindent \textsc{Murotani} Takahiro is assistant professor in anabelian arithmetic geometry at the \emph{Kyoto Institute of Technology}, Japan \textbullet Email: murotani@kit.ac.jp
	
	\medskip
	
	\noindent \textsc{Yamaguchi} Naganori is postdoctoral researcher in anabelian arithmetic geometry at the \emph{Institute of Science Tokyo}, Japan \textbullet Email: yamaguchi.n.ac@m.titech.ac.jp \textbullet \url{https://n-yamaguchi-0729.github.io/homepage-en}
	
}	
\end{document}